\documentclass[11pt]{article}
\usepackage{amssymb,amsfonts,amsmath,curves}
\usepackage{epsfig}
\usepackage{graphicx}
\usepackage{color}
\usepackage{anysize}

\begin{document}

\newtheorem{defin}{Definition}[section]
\newtheorem{Prop}{Proposition}
\newtheorem{theorem}{Theorem}[section]
\newtheorem{ml}{Main Lemma}
\newtheorem{con}{Conjecture}
\newtheorem{cond}{Condition}
\newtheorem{conj}{Conjecture}
\newtheorem{prop}[theorem]{Proposition}
\newtheorem{lem}{Lemma}[section]
\newtheorem{rmk}[theorem]{Remark}
\newtheorem{cor}{Corollary}[section]
\renewcommand{\theequation}{\thesection .\arabic{equation}}

\newcommand{\beq}{\begin{equation}}
\newcommand{\eeq}{\end{equation}}
\newcommand{\beqn}{\begin{eqnarray}}
\newcommand{\beqnn}{\begin{eqnarray*}}
\newcommand{\eeqn}{\end{eqnarray}}
\newcommand{\eeqnn}{\end{eqnarray*}}
\newcommand{\bprop}{\begin{prop}}
\newcommand{\eprop}{\end{prop}}
\newcommand{\bteo}{\begin{teo}}
\newcommand{\bcor}{\begin{cor}}
\newcommand{\ecor}{\end{cor}}
\newcommand{\bcon}{\begin{con}}
\newcommand{\econ}{\end{con}}
\newcommand{\bcond}{\begin{cond}}
\newcommand{\econd}{\end{cond}}
\newcommand{\bconj}{\begin{conj}}
\newcommand{\econj}{\end{conj}}
\newcommand{\eteo}{\end{teo}}
\newcommand{\brm}{\begin{rmk}}
\newcommand{\erm}{\end{rmk}}
\newcommand{\blem}{\begin{lem}}
\newcommand{\elem}{\end{lem}}
\newcommand{\ben}{\begin{enumerate}}
\newcommand{\een}{\end{enumerate}}
\newcommand{\bei}{\begin{itemize}}
\newcommand{\eei}{\end{itemize}}
\newcommand{\bdf}{\begin{defin}}
\newcommand{\edf}{\end{defin}}
\newcommand{\bpr}{\begin{proof}}
\newcommand{\epr}{\end{proof}}

\newenvironment{proof}{\noindent {\em Proof}.\,\,}{\hspace*{\fill}$\halmos$\medskip}

\newcommand{\halmos}{\rule{1ex}{1.4ex}}
\def \qed {{\hspace*{\fill}$\halmos$\medskip}}

\newcommand{\fr}{\frac}
\newcommand{\Z}{{\mathbb Z}}
\newcommand{\R}{{\mathbb R}}
\newcommand{\E}{{\mathbb E}}
\newcommand{\C}{{\mathbb C}}
\renewcommand{\P}{{\mathbb P}}
\newcommand{\N}{{\mathbb N}}
\newcommand{\var}{{\mathbb V}}
\renewcommand{\S}{{\cal S}}
\newcommand{\T}{{\cal T}}
\newcommand{\W}{{\cal W}}
\newcommand{\X}{{\cal X}}
\newcommand{\Y}{{\cal Y}}
\newcommand{\h}{{\cal H}}
\newcommand{\f}{{\cal F}}

\renewcommand{\a}{\alpha}
\renewcommand{\b}{\beta}
\newcommand{\g}{\gamma}
\newcommand{\G}{\Gamma}
\renewcommand{\L}{\Lambda}
\renewcommand{\l}{\lambda}
\renewcommand{\d}{\delta}
\newcommand{\D}{\Delta}
\newcommand{\e}{\epsilon}
\newcommand{\s}{\sigma}
\newcommand{\B}{{\cal B}}
\renewcommand{\o}{\omega}

\newcommand{\nn}{\nonumber}
\renewcommand{\=}{&=&}
\renewcommand{\>}{&>&}
\newcommand{\<}{&<&}
\renewcommand{\le}{\leq}
\newcommand{\+}{&+&}

%Differentiation
\newcommand{\pa}{\partial}
\newcommand{\ffrac}[2]{{\textstyle\frac{{#1}}{{#2}}}}
\newcommand{\dif}[1]{\ffrac{\partial}{\partial{#1}}}
\newcommand{\diff}[1]{\ffrac{\partial^2}{{\partial{#1}}^2}}
\newcommand{\difif}[2]{\ffrac{\partial^2}{\partial{#1}\partial{#2}}}

\definecolor{Red}{rgb}{1,0,0}
\def\red{\color{Red}}

\title{Level $1$ quenched large deviation principle for random walk
in dynamic random environment}

\author{David Campos\thanks{Facultad de Matem\'{a}ticas, Pontificia Universidad
Cat\'{o}lica de Chile, Vicu\~{n}a Mackenna 4860, Macul, Santiago, Chile}\and
Alexander Drewitz\thanks{
Departement Mathematik,
ETH Z\"urich,
R\"amistrasse 101,
CH-8092 Z\"urich,
Switzerland}\and
Alejandro F. Ram\'\i rez$^{*}$\and
Firas Rassoul-Agha\thanks{Department of Mathematics, University of Utah,
155 South 1400 East, Salt Lake City, UT 84109
%E-mail: firas@math.utah.edu,
%URL: www.math.utah.edu/???firas
}\and
Timo Sepp\"al\"ainen\thanks{ Department of Mathematics,
University of Wisconsin-Madison,
419 Van Vleck Hall,
Madison, WI 53706
%E-mail: seppalai@math.wisc.edu,
%URL: www.math.wisc.edu/???seppalai
}
}

\maketitle
%\tableofcontents

\abstract{ Consider a random walk in a  time-dependent random environment
on the lattice $\mathbb Z^d$. Recently, Rassoul-Agha, Sepp\"al\"ainen and
Yilmaz \cite{rsy11} proved a general large deviation principle under
mild ergodicity assumptions on the random
environment   for such a random walk, establishing first
 level $2$ and $3$ large deviation principles.
Here we present two alternative short proofs of the level $1$ large
deviations under mild ergodicity assumptions on the environment: one
for the continuous time case and another one for the discrete time case.
Both proofs provide the existence, continuity and convexity of the rate function.
Our me\-thods are
based
on the use of the sub-additive ergodic theorem as presented
by Varadhan in \cite{v03}.}
%%%%%%%%%%%%%%%%%%%%%%%%%%%%%%%%%%%%%%%%%%%%%%%%%%%%%%%%%%%%%%%%%%%%%%%%%%%%%%%%%%%%%%%%%%%%%%%%%%%%%%%%%%%%%%%%%%%%%%%%%%%%%%%%%%%%%%%%%%

\medskip
{\footnotesize

}

\section{Introduction}
We consider uniformly elliptic random walks in
time-space random environment both in continuous and discrete
time. We present two alternative short proofs
of the level $1$ quenched large deviation principle under
mild conditions on the environment, based on the use
of the sub-additive ergodic theorem as presented
by Varadhan in \cite{v03}.
Previously,   in the discrete
time case, Rassoul-Agha, Sepp\"al\"ainen and Yilmaz
\cite{rsy11}, proved a level 2 and 3 large deviation
principle, from which the level 1 principle can be derived
via contraction.

Let $\kappa_2>\kappa_1>0$. Denote by  $G:=\{e_1,e_{-1},\ldots, e_d,e_{-d}\}$
the set of unit vectors in $\mathbb Z^d$.
Define $\mathcal Q:=\{v=\{v( e):e\in G\}:
\kappa_1\le\inf_{e\in G}v(e)\le\sup_{e\in G}v(e)\le\kappa_2\}$.
Consider a continuous time Markov process  $\omega:=\{\omega_t:t\ge 0\}$
with state space $\Omega_c:=\mathcal Q^{\mathbb Z^d}$,
so that $\omega_t:=\{\omega_t(x):x\in\mathbb Z^d\}$ with
$\omega_t(x):=\{\omega_t(x,e):e\in G\}\in\mathcal Q$.
We call $\omega$ the {\it continuous time environmental process}.
We assume that for each initial condition $\omega_0$,
the process $\omega$ defines a probability
measure $Q_{\omega_0}^c$ on the
Skorokhod space $D([0,\infty);\Omega_c)$.
Let $\mu$ be an invariant  measure for the
environmental process $\omega$ so that for every
bounded continuous function $f:\Omega_c\to\mathbb R$ and $t\ge 0$ we
have that

$$
\int f(\omega_t)d\mu= \int f(\omega_0)d\mu.
$$
Assume that $\mu$ is also invariant under the action of
space-translations.
Furthermore, we define $Q^c_\mu:=\int Q^c_\omega d\mu$, where
with a slight abuse of notation here $\omega\in\Omega_c$.
For a given trajectory $\omega\in D([0,\infty);\Omega_c)$
consider the process $\{X_t:t\ge 0\}$ defined by
the generator

$$
L_sf(x):=\sum_{e\in G}\omega_s(x,e)(f(x+e)-f(x)),
$$
where $s\ge 0$.
We call this process a {\it continuous time random
walk in a uniformly elliptic time-dependent random environment} and denote for
each $x\in\mathbb Z^d$ by
$P_{x,\omega}^c$ the law on  $D([0,\infty);\mathbb Z^d)$ of this
random walk with initial condition $X_0=x$.
We call $P_{x,\omega}^c$ the {\it quenched law}
starting from $x$ of the random walk.

For $x\in\mathbb R^d$, $|x|_2$, $|x|_1$ and $|x|_\infty$ denote
respectively, their Euclidean, $l_1$ and $l_\infty$-norm.
Also, for $r>0$, we define $B_r(x):=\{y\in\mathbb Z^d: |y-x|_2\le r\}$.
Furthermore, given any topological space $T$, we will denote
by $\mathcal B(T)$ the corresponding Borel sets.

We will also consider a discrete version of this model which we define
as follows.
Let $\kappa>0$  and  $R\subset\mathbb Z^d$ finite. Define $\mathcal P:=\{v=\{v( e):e\in R\}:\inf_{e\in R }v(e)
\ge \kappa, \sum_{e\in R}v(e)=1\}$.
Consider a discrete time Markov process  $\omega:=\{\omega_n:n\ge 0\}$
with state space $\Omega_d:=\mathcal P^{\mathbb Z^d}$,
so that $\omega_n:=\{\omega_n(x):x\in\mathbb Z^d\}$ with
$\omega_n(x):=\{\omega_n(x,e):e\in R\}\in\mathcal P$.
We call $\omega$ the {\it discrete time environmental process}.
Let us denote by $Q^d_\omega$ the
corresponding law of the process defined on the
space $\Omega_d^{\mathbb N}$.
Let $\mu$ be an invariant measure for the
environmental process $\omega$ so that for every
bounded continuous function $f:\Omega_d\to\mathbb R$ and $n\ge 0$ we
have that

$$
\int f(\omega_n)d\mu= \int f(\omega_0)d\mu.
$$
Assume that $\mu$ is also invariant under the action of
space-translations.
Furthermore, we define $Q^d_\mu:=\int Q^d_\omega d\mu$.
Given $\omega\in\Omega_d$ and $x\in\mathbb Z^d$,
consider now the discrete time random walk $\{X_n:n\ge 0\}$
with a law $P^d_{x,\omega}$ on $(\mathbb Z^d)^{\mathbb N}$ defined through
$P^d_{x,\omega}(X_0=x)=1$ and the transition probabilities

$$
P^d_{x,\omega}(X_{n+1}=x+e|X_n=x)=\omega_n(x,e),
$$
for $n\ge 0$ and $e\in R$.
We call this process a {\it discrete time random walk in a
uniformly elliptic 
time-space random environment with jump range $R$} and
 call
$P^d_{x,\omega}$  the {\it quenched law} of
the discrete time random walk starting from $x$.
We  will say that {\it $R$ corresponds to the nearest neighbor case}
if $R=\{e\in\mathbb Z^d: |e|_1=1\}$. 
We say that a subset $A\subset\mathbb Z^d$ is {\it convex}
if there exists a convex subset $V\subset\mathbb R^d$ such
that $A=V\cap\mathbb Z^d$, while we say that $A$ is 
{\it symmetric} if $A=-A$. Throughout, we will assume that
the jump range is $R$ is finite, convex and symmetric or that
it corresponds to the nearest neighbor case.

Throughout we will make the following ergodicity assumption.
Note that we do not demand the environment to be necessarily ergodic
under time shifts.

\medskip

\noindent {\bf Assumption (EC)}. Consider the
continuous time
environmental process $\omega$.
For each $s>0$ and $x\in\mathbb Z^d$   define the
transformation
 $T_{s,x}:D([0,\infty);\Omega_c)\to D([0,\infty);\Omega_c)$
by $(T_{s,x}\omega)_t(y):=\omega_{t+s}(y+x)$.
We say that the environmental process $\omega$ satisfies
{\it assumption (EC)} if
 $\{T_{s,x}:s>0,x\in\mathbb Z^d\}$ is an
ergodic family of transformations acting on the space
$(D([0,\infty);\Omega_c),
\mathcal B(D([0,\infty);\Omega_c)), Q_\mu^c)$. In
other words, the latter means that whenever $A\in\mathcal B(D([0,\infty);\Omega_c))$
is such that $T^{-1}_{s,x}A=A$ for every $s>0$ and
$x\in\mathbb Z^d$, then $Q_\mu^c(A)$ is $0$ or $1$.

\medskip

\noindent {\bf Assumption (ED)}. Consider the
discrete time
environmental process $\omega$.
For  $x\in\mathbb Z^d$   define the
transformation
 $T_{1,x}:D([0,\infty);\Omega_d)\to D([0,\infty);\Omega_d)$
by $(T_{1,x}\omega)_n(y):=\omega_{n+1}(y+x)$.
We say that the environmental process $\omega$ satisfies
{\it assumption (ED)}
if $\{T_{1,x}:x\in R\}$ is an ergodic
family of transformations acting on the space $(\Omega_d^{\mathbb N},
\mathcal B( \Omega_d^{\mathbb N}), Q_\mu^d)$.
In other words, whenever $A\in\mathcal B( \Omega_d^{\mathbb N})$
is such that $T^{-1}_{1,x}A=A$ for every $x\in R,$ then $Q_\mu^d(A)$ is $0$ or $1$.

\medskip

\noindent It is straightforward to check that assumption (ED) is equivalent
to asking that  whenever $A\in \mathcal B(\Omega_d^{\mathbb N})$
is such that $A=T^{-1}_{n,x}A$ for every $x \in R$ and $n\in\mathbb N$ then $Q_\mu^d(A)$ is $0$ or $1$.

 In this paper we present a level $1$ quenched large deviation
principle for both
 the con\-tinuous  and the discrete time random walk in
time-space random environment.
It should be noted that the discrete time
version of our result
can be derived via a contraction principle
 from results that
have been obtained in Rassoul-Agha, Sepp\"al\"ainen and Yilmaz \cite{rsy11} establishing
level $2$ and $3$ large
deviations, for discrete time random walks
on time-space random environments and potentials.
There, the authors also derive
 variational expressions for the rate functions.
Nevertheless, the proofs we present here of both
Theorem \ref{one} and \ref{two}, are short and
direct.

\begin{theorem}
\label{one}  Consider a continuous time random walk $\{X_t:t\ge 0\}$
in a uniformly elliptic time-dependent environment $\omega$
satisfying assumption (EC).
Then, there exists a convex continuous rate function $I_c(x):\mathbb R^d\to
[0,\infty)$ such that the following are satisfied.

\begin{itemize}

\item[(i)] For every open set $G\subset\mathbb R^d$ we have that
$Q_\mu^c$-a.s.
$$
\liminf_{t\to\infty}\frac{1}{t}\log P_{0,\omega}^c\left(\frac{X_t}{t}\in G\right)
\ge-\inf_{x\in G}I_c(x).
$$

\item[(ii)] For every closed set $C\subset\mathbb R^d$ we have that
$Q_\mu^c$-a.s.
$$
\limsup_{t\to\infty}\frac{1}{t}\log P_{0,\omega}^c\left(\frac{X_t}{t}\in C\right)
\le-\inf_{x\in C}I_c(x).
$$

\end{itemize}

\end{theorem}
\medskip

\noindent To state the discrete time version of Theorem \ref{one},
we need to introduce some notation.
Let $R_0:=\{0\}\subset\mathbb Z^d$, $R_1:=R$ and for $n\ge 1$ define

$$
R_{n+1}:=\{y\in\mathbb Z^d: y=x+e\ {\rm for}\ {\rm some}\
x\in R_n\ {\rm and}\ e\in R\},
$$
and $U_n:=R_n/n$.
Note that $R_n$ is the set of sites that a random walk with jump range
 $R$ visits with positive probability at time $n$.
We then define $U$ as the set of limit points of the sequence
of sets $\{U_n:n\ge 1\}$, so that

\begin{equation}
\label{def-u}
U:=\{x\in\mathbb R^d: x=\lim_{n\to\infty}x_n\ {\rm for}\ {\rm some}\
 {\rm sequence}\ x_n\in U_n\}.
\end{equation}

\begin{theorem}
\label{two}  Consider a discrete time random walk $\{X_n:n\ge 0\}$
in a uniformly elliptic time-dependent environment $\omega$
satisfying assumption (ED) with jump range $R$.
Assume that either $(i)$
 $R$ is finite, convex, symmetric and there is a neighborhood of
$0$ which belongs to the convex hull of $R$; $(ii)$ or that $R$ corresponds to the nearest neighbor case.
Consider $U$ defined in (\ref{def-u}). Then
$U$ equals the convex hull of $R$ and there exists a convex  rate function $I_d(x):\mathbb R^d\to
[0,\infty]$ such that $I_d(x)\le |\log\kappa |$ for $x\in U$,
$I_d(x)=\infty$ for $x\notin U$,
$I$ is continuous for every $x\in U^o$
and the following are satisfied.

\begin{itemize}

\item[(i)] For every open set $G\subset\mathbb R^d$ we have that
$Q_\mu^d$-a.s.
$$
\liminf_{n\to\infty}\frac{1}{n}\log P_{0,\omega}^d\left(\frac{X_n}{n}\in G\right)
\ge-\inf_{x\in G}I_d(x).
$$

\item[(ii)] For every closed set $C\subset\mathbb R^d$ we have that
$Q_\mu^d$-a.s.
$$
\limsup_{n\to\infty}\frac{1}{n}\log P_{0,\omega}^d\left(\frac{X_n}{n}\in C\right)
\le-\inf_{x\in C}I_d(x).
$$

\end{itemize}

\end{theorem}

\noindent
Both quenched and annealed large deviations for discrete time random walks on random environments which do not
depend on time, have been thoroughly studied in the case
in which $d=1$ (see the reviews of Sznitman  \cite{s04} and Zeitouni \cite{z06}
for both the one-dimensional and multi-dimensional cases). The first quenched
multidimensional result was obtained by Zerner in \cite{z98}
under the so called plain nestling condition, concerning
the law of the su\-pport of the quenched drift (see also \cite{z06}
and \cite{s04}). In \cite{v03}, Varadhan established both a
general quenched and annealed large deviation principle for discrete
time random walks in static random environments via the use
of the subadditive ergodic theorem. In the quenched case,
he assumed uniform ellipticity and the ergodicity assumption (ED).
Subsequently, in his Ph.D. thesis \cite{r06}, Rosenbluth
extended the quenched result of Varadhan under a condition weaker
than uniform ellipticity, along with a variational formula
for the rate function (see also Yilmaz \cite{y08,y09,y09-2}). The method of Varadhan
based on the subadditive ergodic theorem and of Rosenbluth  \cite{r06}, Yilmaz
\cite{y08}
and Rassoul-Agha, Sep\"al\"ainen, Yilmaz \cite{rsy11}, are closely related to
the use of the subadditive ergodic theorem in the
context of non-linear stochastic homogenization (see for example
the paper of dal Maso, Modica \cite{dmm86}).
Closer and more recent examples of stochastic homogenization
for the Hamilton-Jacobi-Bellman equation with static
Hamiltonians via
the subadditive ergodic theorem  are the work of
 Rezakhanlou and Tarver   \cite{rt00} and
of Souganidis \cite{so99} and
in the context of the totally asymmetric simple $K$-exclusion processes and
growth processes the works of Sepp\"al\"ainen in \cite{s99} and Rezakhanlou in
\cite{r02}. Stochastic homogenization for the Hamilton-Jacobi-Bellman
equation with respect to time-space shifts was treated by Kosygina and
Varadhan in \cite{kv08} using change of measure techniques giving
variational expressions for the effective Hamiltonian.

A particular case of
Theorem \ref{one} is the case of a random walk which has
a drift in a given direction on occupied sites and
in another given direction on unoccupied sites,
where the environment is generated by an
attractive spin-flip particle system or a simple
exclusion process (see Avena, den Hollander and Redig \cite{adhr10}
for the case of a one-dimensional attractive spin-flip dynamics, and
also \cite{adhr11,adsv11,dhdss11}).
This case is also included
in the results presented in \cite{rsy11}.
Another particular case of  Theorem \ref{one} is
a continuous time random walk in a static random environment
with a law which is ergodic under spatial translations:
two of these cases are the Bouchaud trap random walk
with bounded jump rates (see for example \cite{bc06}) and the continuous time random
conductances model (see for
example \cite{dfgw89}). Our proof
would also apply to the polymer measure defined by
a continuous time random walk in
time-dependent random environment and bounded random potential
(see \cite{rsy11}).
Note that Theorem \ref{two} does  include the
classical nearest neighbor case
(a nearest neighbor case example is the random walk on a time-space i.i.d.
environment studied by Yilmaz \cite{y09}).

Our proofs are obtained by directly establishing the level $1$ large deviation
principle and is based on the
sub-additive ergodic theorem as used by Varadhan in \cite{v03}.
Let us note, that in \cite{v03}, Varadhan applies sub-additivity
directly to the logarithm of a smoothed up version of the inverse of the
transition probabilities of the random walk, as opposed to the earlier
approach of Zerner \cite{z98} (see also Sznitman \cite{s98}), where sub-additivity is applied to
a generalized Laplace transform of the hitting times of sites
of the random walk forcing to assume the so called nestling property
on the random walk.
While our methods do not give any explicit information about
the rate function, besides its convexity and continuity, the proofs are short and
simple.

We do not know how to define a smoothed up
version of the transition probabilities as is done by Varadhan
in \cite{v03}. We therefore have to prove directly
an equicontinuity estimate for the
transition probabilities of the random walk,
which is  the main difficulty
in the proofs of Theorems \ref{one} and \ref{two}.
 In the
case of Theorem \ref{one} we follow the method presented in \cite{dgrs12}: we first
express
 the transition probabilities
of the walk in terms of those of a simple symmetric
random walk through a Radon-Nykodym derivative,  then
through the use of Chapman-Kolmogorov equation  we
rely on standard large deviation estimates for the continuous time simple
symmetric random walk.

 In  section \ref{s-two} we
present the proof of  Theorem \ref{one} using
the methods developed in  \cite{dgrs12}.
In section \ref{s-three} we  continue with the proof
of Theorem \ref{two} in the case in which the
jump range of the walk $R$ is convex, symmetric
and a neighborhood of $0$ is contained in its convex hull.
In section \ref{s-four} we  prove
 Theorem \ref{two} for the discrete time nearest neighbor case.
Throughout the rest of the paper we will use the notations $c,C,C',C''$
to refer to different positive constants.

\section{Proof of Theorem \ref{one}}
\label{s-two}
For each $s\ge 0$, let $\theta_s:D([0,\infty);\Omega_c) \to D([0,\infty);\Omega_c)$
denote the canonical
time shift.
As in \cite{dgrs12}, we first define
 for each $0\le s<t$ and $x,y\in\mathbb Z^d$ the quantities

$$
e(s,t,x,y):=P_{x,\theta_s\omega}^c\left(X_{t-s}=y\right),
$$
and

$$
a_c(s,t,x,y):=-\log e(s,t,x,y),
$$
where the subscript $c$ in $a_c$ is introduced to distinguish this
quantity from the corresponding discrete time one.
Note that these functions still depend on the realization of $\omega$.
We call $a_c(s,t,x,y)$ the point to point passage function from $x$ to $y$
between times $s$ and $t$. Due to the fact that
we are considering a continuous time random walk, here we do not need
to smooth out the point to point passage functions (see \cite{v03}).
Nevertheless, there is an equicontinuity issue that should be resolved.
 Theorem \ref{one}  will follow directly
from the following shape theorem. A version of this shape theorem
for a random walk in random potential has been established
as Theorem 4.1 in
\cite{dgrs12} (see also Theorem 2.5 of Chapter 5 of Sznitman \cite{s98}).

\begin{theorem}
\label{three}
{\bf [Shape theorem]} There exists a deterministic convex function
$I_c:\mathbb R^d\to [0,\infty)$ such that $Q_\mu^c-a.s.$, for
any compact set $K\subset \mathbb R^d$

\begin{equation}
\label{unif}
\lim_{t\to\infty}\sup_{y\in tK\cap{\mathbb Z^d}}\left|
t^{-1}a_c(0,t,0,y)-I_c\left(\frac{y}{t}\right)\right|=0.
\end{equation}
Furthermore, for any $M>0$, we can find a compact $K\subset\mathbb R^d$
such that $Q_\mu^c-a.s.$

\begin{equation}
\label{comp}
\limsup_{t\to\infty}\frac{1}{t}\log P_{0,\omega}^c\left(\frac{X_t}{t}\notin K\right)\le-M.
\end{equation}

\end{theorem}

\medskip

\noindent
 Let us first see how to derive Theorem \ref{one} from
Theorem \ref{three}. We will first prove the upper bound of part $(ii)$
of Theorem \ref{one}. By (\ref{comp}) of Theorem \ref{three}, we
know that we can choose a compact set $K\subset\mathbb R^d$ such that

$$
\limsup_{t\to\infty}\frac{1}{t}\log P_{0,\omega}^c\left(\frac{X_t}{t}\notin
K\right)<-\inf_{x\in C}I_c(x),
$$
where $C$ is a closed set. It is therefore enough to prove that

$$
\limsup_{t\to\infty}\frac{1}{t}\log P_{0,\omega}^c\left(\frac{X_t}{t}\in
C\cap K\right)\le -\inf_{x\in C}I_c(x).
$$
Now,

\begin{eqnarray*}
&\limsup_{t\to\infty}\frac{1}{t}\log P_{0,\omega}^c\left(\frac{X_t}{t}\in
C\cap K\right)\le \limsup_{t\to\infty}\frac{1}{t}
\sup_{y\in (tC\cap tK)\cap\mathbb Z^d}
 \log P_{0,\omega}^c\left(X_t=y\right)\\
&=
\limsup_{t\to\infty}\frac{1}{t}
 \log P_{0,\omega}^c\left(X_t=y_t\right),
\end{eqnarray*}
where $y_t\in (tC\cap tK)\cap\mathbb Z^d$, is a point that
maximizes $P_{0,\omega}^c(X_t=\cdot)$. Now, by compactness, there is a
subsequence $t_n\to\infty$ such that

$$
\lim_{n\to\infty}\frac{y_{t_n}}{t_n}=:x^*\in C\cap K,
$$
and $\limsup_{t\to\infty}\frac{1}{t}
 \log P_{0,\omega}^c\left(X_t=y_t\right)=
\limsup_{n\to\infty}\frac{1}{t_n}
 \log P_{0,\omega}^c\left(X_{t_n}=y_{t_n}\right)$.
Thus, by the continuity of $I_c$ and by (\ref{unif}) we see that

$$
\limsup_{t\to\infty}\frac{1}{t}\log P_{0,\omega}^c\left(\frac{X_t}{t}\in
C\cap K\right)\le -I_c(x^*)\le -\inf_{x\in C}I_c(x).
$$
To prove the lower bound, part $(i)$ of Theorem \ref{one},
note that by (\ref{unif}) we have that

$$
\liminf_{t\to\infty}\frac{1}{t}
 \log P_{0,\omega}^c\left(\frac{X_t}{t}\in G\right)
\ge
\liminf_{t\to\infty}\frac{1}{t}
\sup_{y\in (tG)\cap\mathbb Z^d}
 \log P_{0,\omega}^c\left(X_t=y\right)\ge -\inf_{x\in G}I_c(x).
$$

\noindent
Let us now continue with the proof of Theorem \ref{three}.
 Display (\ref{comp}) of Theorem \ref{three}
follows from standard large deviation estimates for
the process $\{N_t:t\ge 0\}$, where $N_t$ is the total number
of jumps up to time $t$ of the random walk $\{X_t:t\ge 0\}$,
which can be coupled with a Poisson process of parameter
$2d\kappa_2$.
 To prove the first statement (\ref{unif}) of Theorem \ref{three}
we first observe that for every $0\le t_1<t_2<t_3$ and $x_1,x_2,x_3
\in\mathbb Z^d$ one has that $Q_\mu^c$-a.s.

\begin{equation}
\label{subad}
a_c(t_1,t_3,x_1,x_3)\le a_c(t_1,t_2,x_1,x_2)+a_c(t_2,t_3,x_2,x_3).
\end{equation}
We will also need to obtain bounds
on the point to point passage functions which
will be eventually used to prove some crucial
equicontinuity estimates. To prove these
bounds, we first state Lemma 4.2 of \cite{dgrs12},
which is a large deviation estimate for the simple symmetric
random walk.

\begin{lem}
\label{L:rwldp}
Let $X$  be  a
simple symmetric random walk on $\Z^d$
with jump rate $\kappa$
and starting point
$X(0)=0$.
For each $x\in\mathbb Z^d$ and $t>0$ let $p(t,0,x)$ be the probability
that this random walk is at position $x$ at time $t$ starting from
$0$. Then for every $t>0$ and $x\in\mathbb Z^d$, we have
\beq\label{rwldp}
p(t,0,x)= \frac{e^{-J(\frac{x}{t})\,t}}{(2\pi t)^{\frac{d}{2}}\Pi_{i=1}^d
\big(\frac{x_i^2}{t^2}+\frac{\kappa^2}{d^2}\big)^{1/4}} \left(1+o(1)\right),
\eeq
where
$$
J(x) := \sum_{i=1}^d \frac{\kappa}{d} j\Big(\frac{d x_i}{\kappa}\Big) \qquad
\mbox{with} \quad  j(y) := y\sinh^{-1} y -\sqrt{y^2+1}+1,
$$
and the error term $o(1)$ tends to zero as $t\to\infty$ uniformly in $x\in tK\cap
\Z^d$, for any compact $K\subset \R^d$.
Furthermore the function $j$ is increasing with $|y|$ and $j\ge 0$.
\end{lem}
\medskip
\noindent We will need the following estimates
for the transition probabilities.

\medskip

\begin{lem}
\label{estimate} Consider
the transition probabilities of a random walk on a
uniformly elliptic time-dependent environment.
The following hold $Q_\mu^c$-a.s.

\begin{itemize}

\item[(i)] Let $C_3>0$. There exists a $t_0>0$ and
constants $C_1, C_1'$ and $C_2$ such that
for $\epsilon>0$ small enough and every  $t\ge t_0$,  $y,z\in\mathbb Z^d$
such that $|y-z|_2\le \epsilon t+\frac{tC_3}{|\log\epsilon|}$
 we have that

$$
C_1 e^{-C'_1t\frac{1}{|\log\epsilon|^{1/2}}}
p(\epsilon t,z,y)\le
e(t(1-\epsilon),t,z,y)\le C_2 e^{C_2t\frac{1}{|\log\epsilon|^{1/2}}}
p(\epsilon t,z,y).
$$
\item[(ii)] Let $r>0$. There exists a $t_0>0$ and a constant $C>0$
such that
for each $t\ge t_0$ and $x\in B_{tr}(0)$ one has that

$$
e(0,t,0,x)\ge e^{-Ct}p(t,0,x).
$$

\item[(iii)] There is a function $\alpha: (0,\infty)\times [0,\infty)\to (0,\infty)$ such
that for each $x,y \in \Z^d$ and $t>s \geq 0$ one has that
\beq
\label{le3}
e(s,t,x,y)\ge \alpha(t-s,|x-y|_1)>0.
\eeq
\end{itemize}

\end{lem}
\begin{proof} {\it Part (i)}.
Note that

\begin{equation}
\label{b1}
e(t(1-\epsilon),t,z,y)=E_{z,t(1-\epsilon)}\left[
e^{\int_{t(1-\epsilon)}^t\log(2d\omega_s(Y_{s_-},Y_s-Y_{s_-}))dN_s
-\int_{t(1-\epsilon)}^t (\omega_s(Y_s,G)-1)ds
}1_{Y_t}(y)\right],
\end{equation}
where $E_{z,s}$ is the expectation with respect to the law of
a continuous time simple symmetric random walk $\{Y_t:t\ge 0\}$ of
jump rate $1$ starting from $z$ at time $s$, $N_t$ is
the number of jumps up to time $t$ of the walk, while for each $x\in\mathbb Z^d$
and $s>0$,
$\omega_s(x,G):=\sum_e\omega_s(x,e)$ is the total jump rate
at site $x$ and time $s$ (see for example Proposition 2.6 in Appendix 1 of
Kipnis-Landim
\cite{kl99}).
Using the fact that the jump rates are bounded from above
and from below, it is clear that
there is a constant $C>0$ such that
$$
e^{\int_{t(1-\epsilon)}^t\log(2d\omega_s(Y_{s_-},Y_s-Y_{s_-}))dN_s
-\int_{t(1-\epsilon)}^t (\omega_s(Y_s,G)-1)ds
}\le
e^{C(N_t-N_{t(1-\epsilon)})+C\epsilon t}.
$$
Substituting this bound in (\ref{b1}), we see that

\begin{equation}
\label{ub}
e(t(1-\epsilon),t,z,y)\le
e^{C\epsilon t}E\left[
e^{CN_{\epsilon t}}p_{N_{\epsilon t}}(z,y)\right],
\end{equation}
where now $E$ is the expectation with respect
to a  Poisson process $\{N_t:t\ge 0\}$
of rate $1$ and
 $p_n$ is the $n$-step transition probability of a discrete
time simple symmetric random walk. Let now
$R_\epsilon:=\frac{1}{\epsilon|\log\epsilon|^{1/2}}$. Note that

\begin{eqnarray*}
& E\left[
e^{CN_{\epsilon t}}p_{N_{\epsilon t}}(z,y)\right]\le
e^{C R_\epsilon t\epsilon  }p(\epsilon t,z,y)+E[e^{N_{\epsilon t}C}, N_{\epsilon
    t}> R_\epsilon t\epsilon]\\
&\le
e^{C  R_\epsilon t\epsilon }p(\epsilon t,z,y)+E[e^{2N_{\epsilon t}C}]^{1/2}
P( N_{\epsilon t}> R_\epsilon t\epsilon )^{1/2}.
\end{eqnarray*}
Now, using the exponential Chebychev inequality with parameter
$\log  R_\epsilon$, we get

\begin{equation}
\label{poisson}
P( N_{\epsilon t}> R_\epsilon\epsilon t)\le e^{-\epsilon t ( R_\epsilon\log  R_\epsilon-(R_\epsilon-1))}
\end{equation}
and we compute
$E[e^{2N_{\epsilon t}C}]=e^{\epsilon t(e^{2C}-1)}$. Hence,

\begin{equation}
\label{d6}
E\left[
e^{CN_{\epsilon t}}p_{N_{\epsilon t}}(z,y)\right]\le
e^{C R_\epsilon t\epsilon }p(\epsilon t,z,y)+e^{\epsilon\frac{t}{2}(e^{2C}-1)}e^{-\epsilon
\frac{t}{2} ( R_\epsilon\log  R_\epsilon-( R_\epsilon-1))}.
\end{equation}
Now, by Lemma \ref{L:rwldp} we know that $j(y)$ is increasing
with $|y|$, so that

$$\sup_{y,z:|y-z|_2\le \epsilon t+\frac{C_3t}{|\log\epsilon|}}\epsilon tj\left(\frac{|z-y|}{\epsilon t}\right)
\le \epsilon tj\left(\frac{C_3}{\epsilon|\log\epsilon| }+1\right)\le
t\left(\frac{C_3}{|\log\epsilon|}+\epsilon\right)\log\left(3+\frac{2C_3}{\epsilon|\log\epsilon|}\right)
$$
for $t\ge 1$. Hence, again by Lemma \ref{L:rwldp} with
$\kappa=1$, we see that for any constant $c>0$ we can
choose $\epsilon$ small enough such that

\begin{equation}
\label{d3}
\lim_{t\to\infty}\frac{e^{\epsilon \frac{t}{2}(e^{2C}-1)}
e^{-\epsilon t c(R_\epsilon\log R_\epsilon-(R_\epsilon-1))}}{\inf_{y,z}p(\epsilon t,z,y)}=0,
\end{equation}
where the infimum is taken over $y,z$ as in the previous display. Applying (\ref{d3})  with
$c=1/2$, we see that the second term of the right-hand side
of inequality (\ref{d6}), after taking the supremum
over $y,z$ such that $|y-z|_2\le \epsilon t+\frac{C_3t}{|\log\epsilon|}$,
 is negligible with respect to the first one.
Hence, for $\epsilon$ small enough, there is a constant $C$ and a $t_0>0$ such
that for
$y,z$ such that $|y-z|_2\le \epsilon t+\frac{C_3t}{|\log\epsilon|}$
and $t\ge t_0$ one has
$$
e(t(1-\epsilon),t,z,y)
\le
Ce^{(R_\epsilon+1)Ct\epsilon }p(\epsilon t,z,y).
$$
Similarly, using the fact that the jump rates are bounded from
above and from below it can be shown that for
$y,z$ such that $|y-z|_2\le \epsilon t+\frac{C_3t}{|\log\epsilon|}$
and $t$ large enough
\begin{eqnarray*}
&e(t(1-\epsilon),t,z,y)
\ge
e^{-C'\epsilon t}
E[e^{-C'N_{\epsilon t}}p_{N_{\epsilon t}}(z,y)1_{N_{\epsilon t}\le
    R_\epsilon\epsilon t}]\\
&\ge
e^{-(R_\epsilon+1)\epsilon t C'}E[p_{N_{\epsilon t}}(z,y)1_{N_{\epsilon t}\le R_\epsilon\epsilon t}]
\ge
e^{-(R_\epsilon+1)\epsilon t C'}\left(
p(\epsilon t,z,y)-
P(N_{\epsilon t}> R_\epsilon\epsilon t)\right)\\
&\ge
C''e^{-(R_\epsilon+1)\epsilon t C'}p(\epsilon t,z,y),
\end{eqnarray*}
where we have used (\ref{poisson}) and (\ref{d3}) with $c=1$.

\noindent {\it Part (ii)}. The proof of part $(ii)$ is
analogous to the proof of the lower bound of part $(i)$.

\noindent {\it Part (iii)}. By the same argument as the last part of the proof of
part $(i)$, there is a constant $C'>0$ such that $$e(s,t,x,y) \geq
e^{-C'(t-s)}E[e^{-C'N_{t-s}}p_{N_{t-s}}(x,y),N_{t-s}=|x-y|_1]$$
But $P(N_{t-s}=|x-y|_1)>0$ (there is, with positive probability, a trajectory from
$0$ to $x$ such that $N_{t-s}=|x-y|_1$). Thus,
\begin{eqnarray*}
e(s,t,x,y) &\geq& e^{-C'(t-s)-C'|x-y|_1}p_{|x-y|_1}(x,y)P(N_{t-s}=|x-y|_1)\\
&\ge &e^{-C'(t-s)-C'|x-y|_1}\frac{1}{(2d)^{|x-y|_1}}P(N_{t-s}=|x-y|_1)>0.
\end{eqnarray*}
\end{proof}

\noindent We can now apply Kingman's sub-additive
ergodic theorem (see for example Liggett \cite{l85}),
to prove the following lemma.

\begin{lem} There exists a deterministic function $I_c:\mathbb Q^d\to
[0,\infty)$ such that for every $y\in\mathbb Q^d$, $Q_\mu^c$-a.s.
we have that

\begin{equation}
\label{lQ1}
\lim_{t\to\infty \atop ty\in\mathbb Z^d}\frac{a_c(0,t,0,ty)}{t}=I_c(y).
\end{equation}
\end{lem}
\begin{proof} Assume first that $y\in\Z^d$. Let $q\in\mathbb N$. We will consider for $m > n \geq 1$ the random
variables
$$X_{n,m}(y):=a_c(nq,mq,ny,my).$$
By (\ref{subad}),
 we have
$$
X_{0,m}(y) \leq X_{0,n}(y)+X_{n,m}(y).
$$
By part $(iii)$ of Lemma \ref{estimate}, we see that the
random variables $\{X_{n,m}(y)\}$ are integrable. Hence,
by Kingman's sub-additive ergodic theorem (see Liggett \cite{l85}) we can then
conclude that the limit

\begin{equation}
\label{lQ1.1}
\hat I(q,y,\omega):=\lim_{m \to \infty} \frac{a_c(0,mq,0,my)}{m}
\end{equation}
 exists for $y\in\mathbb Z^d$ and $q\in\mathbb N$. We have to show that it is deterministic. For this reason, let $r>0$,
 $z \in \Z^d$ be arbitrary. It suffices to prove that

$$
\hat I(q,y, \omega) \leq
\hat I(q,y,T_{r,z}\omega)=\lim_{m\to\infty}\frac{a_c(r,mq+r,z,my+z)}{m}.
$$
First, we have that

$$
\frac{a_c(0,mq,0,my)}{m} \leq \frac{a_c(0,r,0,z)}{m}+\frac{a_c(r,mq,z,my)}{m}.
$$
By part $(iii)$ of Lemma \ref{estimate}, the first term of the right-hand side of the last equation tends to
$0$ as $m \to \infty$. Therefore,
\begin{equation}
\label{lQ3}
\hat I(q,y,\omega)=\lim_{m \to \infty}\frac{a_c(0,mq,0,my)}{m} \leq \liminf_{m \to
\infty}\frac{a_c(r,mq,z,my)}{m}.
\end{equation}
On the other hand, for $u\in\mathbb N$ such that $m>u>r$ we have that
\begin{eqnarray*}
\frac{a_c(r,mq,z,my)}{m}
&\leq& \frac{a_c(r,(m-u)q+r,z,(m-u)y+z)}{m}\\
&+&\frac{a_c((m-u)q+r,mq,(m-u)y+z,my)}{m}.\\
\end{eqnarray*}
 Again, by
 part $(iii)$ of Lemma \ref{estimate}, the last term tends to $0$
as $m \to \infty$. Therefore
\begin{equation}
\label{lQ4}
\liminf_{m \to \infty} \frac{a_c(r,mq,z,my)}{m}\leq \lim_{m
\to \infty}\frac{a_c(r,(m-u)q+r,z,(m-u)y+z)}{m}=\hat I(q,y,T_{r,z}\omega).
\end{equation}
Hence  $\hat I(q,y,\omega) \leq \hat I(q,y,T_{r,z}\omega)$. Since $r>0$ and $z \in \Z^d$
are arbitrary, $\hat I(q,y)$ is shift-invariant under each transformation $T_{r,z}$.
By assumption (EC), $\hat I(q,y)$ is $Q_\mu^c$-a.s equal to a constant for
each $y$.
Now, if $y\in\mathbb Q^d$, choose the smallest $q\in\mathbb N$
such that $qy\in\mathbb Z^d$. Then by  (\ref{lQ1.1}), we conclude that

\begin{equation} \label{eq:IcDef}
\lim_{m \to \infty} \frac{a_c(0,mq,0,mqy)}{mq}=\frac{1}{q}\hat I(q,qy,\omega)
=:I_c(y),
\end{equation}
exists (and is well-defined) and is $Q_\mu^c$-a.s. equal to a constant.
\end{proof}

\noindent We now need to extend the definition of the function $I_c(x)$
for all $x\in\mathbb R^d$ and prove the uniform convergence
in (\ref{unif}). To do this, we will prove that for each compact $K$
there is a $t_0>0$ such that the
family of functions $\{t^{-1}a_c(0,t,0,ty): t\ge t_0\}$ defined on $K$
is equicontinuous.
We can now proceed to  the main step of the
proof of Theorem \ref{three}.

\begin{lem}
\label{L:modcont}
Let $K$ be any compact subset of $\R^d$. There exist deterministic $\phi_K:
(0,\infty)\to(0,\infty)$ with $\lim_{r\downarrow 0}\phi_K(r)=0$, and  $t_0>0$
such that
for any $\epsilon>0$ and $t\ge t_0$, $Q_{\mu}^c$-a.s., we have
\beq\label{modcont1}
\sup_{x,y\in tK \cap \Z^d \atop  |x-y|_2 \leq \epsilon t}
t^{-1}|a_c(0,t,0,x)-a_c(0,t,0,y)| \leq \phi_K(\epsilon).
\eeq
\end{lem}
\begin{proof} Let us note that for every $\epsilon>0$, $t$ and
$x\in\mathbb Z^d$ one has that

$$
e(0,t,0,x)=\sum_{z\in\mathbb Z^d}e(0,t(1-\epsilon),0,z)e(t(1-\epsilon),t,z,x).
$$
Let $R_K:=\sup\{|x|_2:x\in K\}$ be the maximal distance to $0$
for any point in $K$  and
$r_K=\frac{C_K}{|\log\epsilon|}$, where $C_K$ is a constant that
will be chosen large enough.
From part $(i)$ of Lemma \ref{estimate} and Lemma \ref{L:rwldp}, note that
for $t\ge t_0$ (where $t_0$ is given by part $(i)$ of Lemma \ref{estimate})
\begin{equation}
\label{in5}
e(0,t,0,x)
\le
\sum_{ z\in B_{r_Kt}(x)}
e(0,t(1-\epsilon),0,z)e(t(1-\epsilon),t,z,x)
+
C e^{\frac{1}{|\log\epsilon|^{1/2}} t C-\epsilon
t\frac{1}{d}j\left(d\frac{r_K}{\epsilon}\right)}.
\end{equation}
On the other hand by part $(ii)$ of Lemma \ref{estimate} we have that
for $t\ge t_0$

$$
e(0,t,0,x)\ge e^{-C't-tJ\left(\frac{x}{t}\right)}.
$$
Using the upper bound $J\left(\frac{x}{t}\right)\le dR_K\log (1+2dR_K)$
we see that if

\begin{equation}
\label{in10}
\epsilon\frac{1}{d} j\left(d\frac{r_K}{\epsilon}\right)> C+C'+
dR_K\log\left(1+2dR_K\right),
\end{equation}
the second term of (\ref{in5}) is negligible. But (\ref{in10})
is satisfied for $C_K>2(C+C'+dR_K\log (1+2dR_K))$
and $\epsilon>0$ small enough.
Hence, it is enough to prove that, $Q_\mu^c$-a.s. we have that

\beq\label{modcont2}
\sup_{x,y\in tK \cap \Z^d \atop |x-y|_2 \leq \epsilon t}
\sup_{z\in B_{r_Kt}(x)}
\frac{e(t(1-\epsilon),t,z,x)}{e(t(1-\epsilon),t,z,y)}\le
 Ce^{t\phi_K(\epsilon)}.
\eeq
To this end, by Lemmas \ref{L:rwldp}  and \ref{estimate}

\begin{equation}
\label{eq1}
\frac{e(t(1-\epsilon),t,z,x)}{e(t(1-\epsilon),t,z,y)}
\le Ce^{2tC\frac{1}{|\log\epsilon|^{1/2}}}e^{-\epsilon t
\left(J\left(\frac{x-z}{\epsilon t}\right)-J\left(\frac{y-z}{\epsilon t}
\right)\right)}.
\end{equation}
But,

\begin{align}
\nonumber
&J\left(\frac{z-x}{t\epsilon}\right)-J\left(\frac{z-y}{t\epsilon}\right)
=\sum_{i=1}^d\frac{1}{d}\left[j\left(d
\frac{z_i-x_i}{t\epsilon}
\right)-j\left(d\frac{z_i-y_i}{t\epsilon}
\right)\right]\\
\nonumber
&\le
\sum_{i=1}^d\left|\frac{1}{d}\int_{d
\frac{z_i-x_i}{t\epsilon}}^{d\frac{z_i-y_i}{t\epsilon}}
\log\left(1+2|u|\right)du\right|
\le d\log\left(1+\frac{2dC_K}{\epsilon |\log\epsilon|} \right).
\end{align}
Substituting this estimate back into (\ref{eq1}) we obtain
(\ref{modcont2})
with $\phi_K(\epsilon)=C\frac{1}{|\log\epsilon |^{1/2}}$.
\end{proof}
Using this lemma, we can extend $I_c$ to a continuous function on $\R^d.$
It remains to show the convexity of $I_c.$ For this purpose, let $\lambda \in (0,1),$
$x, y \in \R^d$ and let $(\lambda_n) \subset (0,1) \cap \mathbb Q,$ $(x_n), (y_n) \subset \mathbb Q^d$
such that $\lambda_n \to \lambda,$ $x_n \to x,$ and $y_n \to y.$
In addition let $r_n \in \N$ be such that
$
r_n(\lambda_n x_n + (1-\lambda_n)y_n),$ $\lambda_n mr_n,$ and $\lambda_n mr_n x_n,$
are contained in $\Z^d.$
Then for any $n \in \N$ one has
\begin{align*}
 I_c(\lambda_n x_n + (1-\lambda_n)y_n)
&= \lim_{m \to \infty} \frac{a_c(0,mr_n,0,mr_n(\lambda_n x_n + (1-\lambda_n)y_n))}{mr_n}\\
&\le \lim_{m \to \infty}  \frac{a_c(0,\lambda_n mr_n,0,\lambda_n mr_n x_n )}{mr_n}\\
&\quad + \lim_{m \to \infty} \frac{a_c(\lambda_n mr_n, mr_n,\lambda_n mr_n x_n, mr_n(\lambda_n x_n+ (1-\lambda_n)y_n)) }{mr_n}. \\
\end{align*}
Now taking $n \to \infty,$
the continuity of $I_c$ yields that the left-hand side converges to $I_c(\lambda x + (1-\lambda)y).$
Taking advantage of the continuity of $I_c$ and \eqref{eq:IcDef},
 the first summand on the right-hand side converges to $\lambda I_c(x)$ a.s., while in combination with the fact that
the transformations $T_{\lambda_n mr_n, \lambda_n mr_n x_n}$ are measure preserving, the second summand
converges in probability to $(1-\lambda) I_c(y);$ from the last fact we deduce a.s. convergence along an appropriate
subsequence and hence the convexity of $I_c.$

\section{Proof of Theorem \ref{two} for the convex case}
\label{s-three}
Here we consider the case in which the jump range $R$ of the walk is convex,
symmetric and a neighborhood of $0$ is contained in the convex hull
of $R$.
Let us call $\pi_{n,m}(x,y)$,
the probability that the discrete time random walk
in time-space random environment
jumps from time $n$ to time $m$ from site $x$ to site $y$.
Define

$$
a_d(n,m,x,y):=-\log \pi_{n,m}(x,y).
$$
As in the continuous time case, we have
the following sub-additivity property
for $n\le p\le m$ and $x,y,z\in\mathbb Z^d$,

\begin{equation}
\label{subd}
a_d(n,m,x,y)\le a_d(n,p,x,z)+a_d(p,m,z,y).
\end{equation}
We first need to define some  concepts
that will be used throughout this section.
An element $(n,z)$ of the set $\mathbb N\times \mathbb Z^d$
will be called a {\it time-space point}. The
time-space points of the form $(1,z)$, with $z\in R$, will be called
{\it steps}. Furthermore,
given two time-space points $(n_1,x^{(1)})$ and $(n_2,x^{(2)})$ a
sequence of steps $(1,z^{(1)}),\ldots, (1,z^{(k)})$, with
$k=n_2-n_1$ will be
called an {\it admissible path from $(n_1,x^{(1)})$ to $(n_2,x^{(2)})$}, if
 $x^{(2)}=x^{(1)}+z^{(1)}+\ldots+ z^{(k)}$ and

\begin{eqnarray}
\nonumber
&\pi_{n_1,n_1+1}(x^{(1)},x^{(1)}+z^{(1)})
\pi_{n_1+1,n_1+2}(x^{(1)}+z^{(1)},x^{(1)}+z^{(1)}+z^{(2)})\times\cdots\\
\label{admiss}
&\cdots \times\pi_{n_2-1,n_2}(x^{(1)}+z^{(1)}+\cdots +z^{(k-1)},x^{(1)}
+z^{(1)}+\cdots +z^{(k)})>0.
\end{eqnarray}
In other words, there is a positive probability
for the time-space random walk $(n,X_n)$ to jump
through the sequence of time-space points $(n_1,x^{(1)}), (n_1+1,x^{(1)}
+z^{(1)}),
\ldots ,(n_2,x^{(2)})=(n_2,x^{(1)}+z^{(1)}+\cdots+z^{(k)})$.
Note that the sequence of steps $(1,z^{(1)}),\ldots, (1,z^{(k)})$, 
is an admissible path if and only if $z^{(j)} \in R$ for all $1 \le j \le k$.
Let us note that by uniform ellipticity asking that
the left-hand side of (\ref{admiss}) be positive
is equivalent to asking that it be larger than or equal
to $\kappa^{n_2-n_1}$. With a slight abuse of notation,
we will adopt the convention that for $u\in\mathbb R$, $[u]$
is the integer closest to $u$ that is between $u$ and $0$.
Furthermore, we introduce for $x\in\mathbb R^d$, the notation
$[x]:=([x_1],\ldots,[x_d])\in\mathbb Z^d$.
 Throughout, given $A\subset \mathbb R^d$ we will call
$A^o$ its interior.

\smallskip

\begin{lem}
\label{ru} Consider a discrete time random walk in a uniformly
elliptic time-dependent environment $\omega$ with
finite, convex and symmetric jump range $R$ such that
a neighborhood of $0$ belongs to its convex hull. Then, $U$ equals the
convex hull of $R$ and for every $n\ge 1$ we have that

\begin{equation}
\label{ru1}
 R_n= (nU)\cap\mathbb Z^d.
\end{equation}
\end{lem}
\begin{proof} It is straightforward to check
that $U$ equals the convex hull of $R$ in $\mathbb R^d$.
On the other hand, note that  if $x\in R_n$, we have that for every $m\in\mathbb N$,
$mx\in R_{nm}$, which implies that $\frac{x}{n}\in U_{nm}$.
This proves  that $R_n\subset (nU)\cap\mathbb Z^d$. Finally, using the
fact that $R$ is convex, we can prove that $(nU)\cap \mathbb Z^d\subset R_n$.
\end{proof}

\medskip

\noindent For each $x\in\mathbb Z^d$ define 
$s(x)$ as the minimum number $n$ of steps
such that there is an admissible path between $(0,0)$
and $(n,x)$. Alternatively, 
 
$$
s(x)=\min\{n\ge 0: x\in R_n\}.
$$
Let us now define a norm in $\mathbb R^d$ which will be a
good approximation for the previous quantity. For each
 $y\in \partial U$ define $||y||=1$. Then, for each $x\in \mathbb R^d$
which is of the form $x=ay$ for some real $a\ge 0$, we define
$||x||=a$. Note that since $U$ is convex, symmetric and
there is a neighborhood of $0$ which belongs to its
interior, this defines a norm in
$\mathbb R^d$ (see for example Theorem 15.2  of Rockafellar \cite{r97}) and that
$x\in U^o$ if and only if $||x||<1$. Furthermore, note
that for every $x\in\mathbb R^d$ we have that

\begin{equation}
\label{norm}
||x||\le s(x)\le ||x||+1.
\end{equation}

\medskip

\begin{lem}
\label{s3lem}
Let $z\in U$ and $x \in U^o$.  Then, for each natural $n$ there exists
an $n_2$ such that 
\begin{equation}
\label{lemd1}
n \leq n_2 \leq n+\frac{9}{1-||x||}+n\frac{||x-z||}{1-||x||}.
\end{equation}
and there is an admissible path between $(n,z)$ and $(n_2,x)$
so that
\begin{equation}
\label{lemd2}
a_d(0,n_2,0,[n_2 x])\leq a_d(0,n,0,[nz])-\log \kappa^{n_2-n}.
\end{equation}
Similarly, for each natural $n$ there exists an $n_1$ such that
\begin{equation}
\label{lemd3}
 n-\frac{9}{1-||x||}-n\frac{||x-z||}{||1-x||}\le n_1 \leq n
\end{equation}
and there is an admissible path between $(n_1,x)$ and $(n,z)$ so that
\begin{equation}
\label{lemd4}
a_d(0,n,0,[n z])\leq a_d(0,n_1,0,[n_1x])-\log \kappa^{n-n_1}
\end{equation}
\end{lem}

\begin{proof} Assume that $n_2\ge n$. It is enough to prove that for $n$ and $n_2$
satisfying (\ref{lemd1}) and (\ref{lemd2}) it is true  that

\begin{equation}
\label{lemd5}
s\left([n_2x]-[nz]\right)
\le n_2-n.
\end{equation}
Now, by (\ref{norm}) and the fact that $||x-[x]||\le 2$ we have that

\begin{eqnarray*}
\nonumber
& 
s\left([n_2x]-[nz]\right)\le
||[n_2x]-[nz]||+1\le ||[n_2x]-[nx]||+||[nx]-[nz]||+1\\
&\le ||(n_2-n)x||+||n(x-z)||+9=(n_2-n)||x||+n||x-z||+9.
\end{eqnarray*}
It follows that to prove (\ref{lemd5}) it is enough to show that

\begin{equation}
\label{lemd6}
(n_2-n)||x||+n||x-z||+9\le n_2-n,
\end{equation}
which is equivalent to 

$$
n_2\ge n+\frac{9}{1-||x||}+n\frac{||x-z||}{1-||x||}.
$$
This proves (\ref{lemd1}). Now assume that $n_1\le n$. We have to show
that

$$
s\left([nz]-[n_1x]\right)
\le n-n_1.
$$
Now,

$$s\left([nz]-[n_1x]\right)\le
||[nz]-[n_1x]||+1\le n||z-x||+(n-n_1)||x||+9.
$$
Hence, it is enough to show that

$$
n||z-x||+(n-n_1)||x||+9\le n-n_1,
$$
which is equivalent to

$$
n_1\le n-\frac{9}{1-||x||}-n\frac{||z-x||}{||1-x||}.
$$
\end{proof}

\smallskip

\noindent We are now ready to prove the following proposition.

\medskip
\begin{prop}
\label{propos} For each $x\in\mathbb R^d$ we have that
$Q_\mu^d$-a.s.
the limit

$$
I(x):=-\lim_{n\to\infty}\frac{1}{n}\log\pi_{0,n}(0,[nx]),
$$
exists, is convex and  deterministic. Furthermore, $I(x)<\infty$
if and only if $x\in U$.
\end{prop}

\begin{proof}
From Lemma \ref{ru}, it follows that for $x\notin U$
it is true for $n\ge 1$, that $nx\notin nU$ and hence from
Lemma \ref{ru} that $nx\notin R_n$
so that $\pi_n(0,[nx])=0$.
Thus,
$I(x)=\infty$.
We divide the rest of the proof in  four steps.
In step $1$
 for each $x\in\mathbb Q^d\cap U^o$ we define
a  function $\tilde I(x)$. In step $2$ we will
show that $\tilde I$ is deterministic
for $x\in\mathbb Q^d\cap U^o$. In step $3$ we will show that
$I(x)$ is well-defined for
$x\in\mathbb Q^d\cap U^o$ and  that $I(x)=\tilde I(x)$
and in step $4$, we extend the definition of $I(x)$
to $x\in U$.

\smallskip
\noindent {\it Step 1.} Here we will define
 for each $x\in\mathbb Q^d\cap U^o$
a  function $\tilde I(x)$.
Given $x\in\mathbb Q^d\cap U^o$, there exist
a $k\in\mathbb N$ and
a $y\in\mathbb Z^d\cap kU^o$ such that
$x=k^{-1}y$. Now, by 
display (\ref{ru1}) of Lemma \ref{ru} we know
that $y\in R_k$. Then, by the convexity of $R$ and the sub-additive
ergodic theorem and (\ref{subd}) we can define $Q^d_\mu$-a.s.

$$
\tilde I(k^{-1}y):=-\lim_{m\to\infty}\frac{1}{mk}\log\pi_{0,mk}(0,my).
$$
This definition is independent of the representation
of $x$. Indeed, assume that $x=k^{-1}y_1=l^{-1}y_2$
for some $k,l\in\mathbb N$, $y_1\in \mathbb Z^d\cap kU^o$
and $y_2\in \mathbb Z^d\cap lU^o$.
Then, passing to subsequences,

\begin{eqnarray*}
\tilde I(k^{-1}y_1)&=&-\lim_{n\to\infty}\frac{1}{nlk}\log\pi_{0,nlk}(0,nly_1)\\
&=&-\lim_{n\to\infty}\frac{1}{nlk}\log\pi_{0,nlk}(0,nky_2)
=\tilde I(l^{-1}y_2).
\end{eqnarray*}

\smallskip
\noindent {\it Step 2.} Here we will show that $\tilde I$ is deterministic in
$\mathbb Q^d\cap U^o$. Let $x \in \mathbb Q^d\cap U^o$.
We know that there exists a $k\in\mathbb N$ and a $y\in\mathbb Z^d
\cap kU^o$ such that $x=k^{-1}y$. Let us now fix $z \in R$. It suffices to prove that
$$
\tilde I(x, \omega) \leq
\tilde I(x,T_{1,z}\omega)=\lim_{m\to\infty}\frac{a_d(1,mk+1,z,my+z)}{mk}.
$$
First, for each $n \in \N$, we have that

$$
\frac{a_d(0,mnk,0,mny)}{mnk} \leq \frac{a_d(0,1,0,z)}{mnk}+\frac{a_d(1,mnk,z,mny)}{mnk}.
$$
 By uniform ellipticity, the first term of the right-hand side of the last inequality tends to
$0$ as $m \to \infty$. Therefore,
\begin{equation}
\label{S21}
\tilde I(x,\omega)=\lim_{m \to \infty}\frac{a_d(0,mnk,0,mny)}{mnk} \leq \liminf_{m \to
\infty}\frac{a_d(1,mnk,z,mny)}{mnk}.
\end{equation}
On the other hand,
\begin{eqnarray}
\nonumber
\frac{a_d(1,mnk,z,mny)}{mnk}
&\leq& \frac{a_d(1,(m-1)nk+1,z,(m-1)ny+z)}{mnk}\\
\label{s21a}
&+&\frac{a_d((m-1)nk+1,mnk,(m-1)ny+z,mny)}{mnk}.
\end{eqnarray}
Let us now assume that
there is an admissible path from $(0,z+(m-1)ny)$ to
$(nk-1,mny)$. This is equivalent to asking
that $z$ satisfies the following condition:

\begin{equation}
\label{cond}
\pi_{0,nk-1}(z+(m-1)ny,mny)>0\qquad {\rm for\ some}\ n\in\mathbb N.
\end{equation}
Then,  by uniform ellipticity, the last term of (\ref{s21a})
tends to $0$
as $m \to \infty$. Therefore, if $z\in R$
satisfies condition (\ref{cond}), by (\ref{S21}) and (\ref{s21a}) we have that

\begin{equation}
\label{s22}
\tilde I(x,\omega)\le \tilde I(x,T_{1,z}\omega).
\end{equation}
Hence, to finish the proof it is enough to show that every $z\in R$
satisfies (\ref{cond}). Now,
$z$ satisfies (\ref{cond}) if and only if
there exists an $n\in\mathbb N$ such that

\begin{equation}
\label{cond2}
z-ny\in R_{nk-1}.
\end{equation}
We will show by contradiction that every $z\in R$ satisfies
(\ref{cond2}). Indeed, assume that for each $n$ it
is true that

$$z-ny\notin R_{nk-1}.$$
Then,

$$\frac{z}{nk-1}-y\frac{n}{nk-1}\notin U_{nk-1}.$$
Therefore,  taking the limit $n \to \infty$, we 
conclude that   $\frac{y}{k} \notin U^o$, which is a contradiction.
This proves that for every $z\in R$  condition (\ref{cond})
is satisfied and hence (\ref{s22}) is also valid.
It follows now by the ergodicity assumption (ED), that
for each $x \in \mathbb Q^d \cap U^o$, $\tilde I(x)$ is $Q^d_\mu$-a.s equal to a constant.

\smallskip
\noindent {\it Step 3.} Here we will show that $I$ is well-defined
in $\mathbb Q^d\cap U^o$ and hence equals $\tilde I$ there.
Let $x\in \mathbb Q^d\cap U^o$. Let $k$ be such that
$kx\in\mathbb Z^d$. Given $n$, choose $m$ so that
$mk\le n < (m+1)k$. Note that
there exists a sequence of increments
$z^{(j)}\in R$, $1\le j\le n-mk$,
such that

$$
[nx]=mkx+z^{(1)}+\cdots +z^{(n-mk)}.
$$
Hence, by sub-additivity and considering that by uniform
ellipticity the  path
$(1,z^{(1)}),\ldots,$  $(1,z^{(n-mk)})$ from $[nx]$ to $mkx$ is
admissible,
  we conclude
that

$$
\frac{a_d(0,n,0,[nx])}{n}\le \frac{a_d(0,mk,0,mkx)}{n}
-\frac{\log \kappa^{n-mk}}{n}.
$$
It follows that

$$
\limsup_{n\to\infty}\frac{a_d(0,n,0,[nx])}{n}\le \tilde I(x).
$$
For the upper bound, first note that similarly
there exists an admissible path of $(m+1)k-n$ steps from $[nx]$
to $(m+1)kx$. Hence,

$$
\frac{a_d(0,(m+1)k,0,(m+1)kx)}{n}\le\frac{a_d(0,n,0,[nx])}{n}
-\frac{\log\kappa^{(m+1)k-n}}{n}.
$$
Taking the limit when $n\to\infty$ we obtain

$$
\liminf_{n\to\infty}\frac{a_d(0,n,0,[nx])}{n}\ge \tilde I(x).
$$

\smallskip
\noindent {\it Step 4.} Here we will show that $I$ is well-defined
in the set $(\mathbb R^d\backslash \mathbb Q^d)\cap U^o$.
 Let $z \in (\mathbb R^d\backslash \mathbb Q^d)\cap U^o$. Pick a rational point $x$ such
that 

\begin{equation}
\label{ff}
\frac{1}{1-||x||}
\le
2\frac{1}{1-||z||}.
\end{equation}
For each $n$,
 from Lemma \ref{s3lem}, we can  find $n_1,n_2$ such that
$n_1 \leq n \leq n_2$,

$$\frac{n_2}{n}\cdot \frac{1}{n_2}a_d(0,n_2,0,[n_2x])\leq \frac{1}{n}a_d(0,n,0,[nz])+b \left(\frac{n_2}{n}-1\right)$$
and

$$\frac{1}{n}a_d(0,n,0,[nz])\leq \frac{n_1}{n}\cdot \frac{1}{n_1}a_d(0,n_1,0,[n_1x])+b\left(1-\frac{n_1}{n}\right),$$
where $b=-\log \kappa \in (0,\infty)$. Take $n \to \infty$. From (\ref{lemd1}) and (\ref{lemd3}) and taking
$C(z)=2\frac{1}{1-||z||}$,
 the limit points of $\frac{n_2}{n}-1$ and $1-\frac{n_1}{n}$ lie in the
 interval $[0,C(z)||x-z||]$
because $x$ satisfies (\ref{ff}). Consequently from the last two inequalities
we see that

\begin{equation}
\label{lcont1}
I(x)\leq \liminf_{n \to \infty}\frac{1}{n} a_d(0,n,0,[nz])+C(z)b||x-z||
\end{equation}
and

\begin{equation}
\label{lcont2}
\limsup_{n \to \infty} \frac{1}{n}a_d(0,n,0,[nz])\leq
I(x)+C(z)b||x-z||.
\end{equation}
Letting $x \to z$, we  conclude that $I$ is well-defined in the set $(\mathbb R^d\backslash \mathbb Q^d)\cap U^o$.

\end{proof}

\smallskip

\noindent We are now in a position to introduce the rate function
of Theorem \ref{two}. We define, for each $x\in U$,

\begin{equation}
\label{defi}
I_d(x):=\left\{
\begin{array}{ll}
I(x) & {\rm for}\ x\in U^o\\
\liminf_{U^o\ni y\to x} I(y) & {\rm for}\ x\in \partial U \\
\infty & {\rm for}\ x\notin U.
\end{array}
\right.
\end{equation}
We will now prove that $I_d$ satisfies the requirements of Theorem \ref{two}.
By uniform ellipticity, it is clear that $I(x)\le|\log\kappa|$
when $x\in U$. From (\ref{lcont1}) and (\ref{lcont2}),
we see that $I$ is continuous in the interior of $R$ (in fact,
Lipschitz continuous in any compact contained in $U^o$).
%Now, by
%Theorem 10.2 of Rockafellar \cite{r97}, we conclude that
%$I$ is upper semi-continuous in $U$.
These observations imply
that $I_d$ defined in (\ref{defi}) is bounded
by $|\log\kappa|$ in $U$, is continuous in $U^o$,
and is lower semi-continuous in $U$.
The convexity of $I_d$ is derived in a manner similar to the continuous time case.
We now prove parts $(i)$ and $(ii)$
of Theorem \ref{two}.

Part $(i)$ of Theorem \ref{two} follows immediately from
the definition of $I_d$ and the fact that for open
sets $G$, $\inf_{x\in G}I(x)=\inf_{x\in G}I_d(x)$.
To prove part $(ii)$ we first consider a compact set $C$ contained in  $U^o$. In this case, we have
\begin{eqnarray*}
&\limsup_{n \to \infty}\frac{1}{n}\log P_{0,\omega}^d \left(\frac{X_n}{n} \in C\right)\leq
\limsup_{n \to \infty}\sup_{x \in C}\frac{1}{n}\log \pi_{0,n}(0,[nx])\\
&=\inf_n \sup_{m \geq n}\sup_{x \in C}\frac{1}{m}\log \pi_{0,m}(0,[mx])
=\inf_{n}\sup_{x \in C}\sup_{m \geq n}\frac{1}{m}\log \pi_{0,m}(0,[mx])\\
&= \inf_n\sup_{x \in C} a_n(x),
\end{eqnarray*}
where we have defined for $x\in U^o$,

$$
a_n(x):=\sup_{m \geq n}\frac{1}{m}\log \pi_{0,m}(0,[mx]).
$$

%Now, by Proposition \ref{propos}, we know  that $\lim_{n\to\infty}a_n(x)= -I(x)$.

\noindent  Hence, the upper bound follows if we can show that, for any given $\epsilon>0$,

$$
\sup_{x \in C}a_n(x)\leq -\inf_{x \in C}I(x)+\epsilon
$$
for large enough $n$. If we assume the opposite, we
can find points $z_m\in C$ which have
a subsequence converging to $z\in C$
and such that along this subsequence one also has that

$$
\frac{1}{m}\log \pi_{0,m}(0,[mz_{m}])>-I(z)+\epsilon.
$$
\noindent
Applying the first part of Lemma \ref{s3lem} gives an index $m_2>m$ such that

$$
\frac{1}{m_2}\log \pi_{0,m_2}(0,[m_2z])\geq
\frac{m}{m_2}(-I(z)+\epsilon)-b\left(1-\frac{m}{m_2}\right).
$$
\noindent Now, since  $\lim_{m\to\infty}\frac{m}{m_2} = 1$
and since by Proposition \ref{propos} $\lim_{m_2\to\infty}\frac{1}{m_2}
\log\pi_{0,m_2}(0,[m_2z])=-I(z)$,
we obtain that $-I(z) \geq -I(z)+\epsilon$, which is a contradiction.

\medskip

\noindent In the general case, let $C \subset U$ be a compact set. Fix $\delta>0$ and let $C_1=\frac{1}{1+\delta}C$. Now $C_1$ is a compact set contained in $U^o$. Pick $\epsilon>0$ small enough so that the closed $\epsilon-$fattening
$C_2=\overline{C_1^{(\epsilon)}}$ is still a compact set contained in $U^o$. Let $n_2=\lfloor(1+\delta)n\rfloor$. Then for large enough $n$, $\frac{x}{n} \in C$ implies $\frac{x}{n_2} \in C_2$. By uniform ellipticity, we have that

\begin{eqnarray*}
&P_{0,\omega}^d \left(\frac{X_n}{n} \in C\right) \kappa^{n_2-n} = \sum_{x \in nC \cap \Z^d}P_{0,\omega}^d(X_n=x)\kappa^{n_2-n}\\
&\leq \sum_{x \in nC \cap \Z^d}P_{0,\omega}^d(X_n=x)\pi_{n,n_2}(x,x)
=\sum_{x \in nC \cap \Z^d}P_{0,\omega}^d(X_n=x,X_{n_2}=x)\\
& \leq \sum_{x \in nC \cap \Z^d}P_{0,\omega}^d(X_{n_2}=x)
 \leq P_{0,\omega}^d \left(\frac{X_{n_2}}{n_2} \in C_2\right),
\end{eqnarray*}

\noindent where the last inequality is satisfied for $n$ large enough. Then,
from the first step of the proof of part $(ii)$ of Theorem \ref{two}
$$\limsup_{n\to \infty}\frac{1}{n}\log P_{0,\omega}^d\left(\frac{X_n}{n} \in C\right) \leq -\inf_{x \in C_2}I(x)+\delta b.$$

\noindent By taking $\epsilon \searrow 0$ and using compactness
and the continuity of $I$

$$\limsup_{n\to \infty}\frac{1}{n}\log P_{0,\omega}^d\left(\frac{X_n}{n} \in C\right) \leq -\inf_{x \in C_1}I(x)+\delta b.$$

\noindent Take $\delta \searrow 0$ along a subsequence $\delta_j$. This takes $C_1$ to C. For each $\delta_j$,
let $z_j \in C_1=C_1(\delta_j)$ satisfy $I(z_j)=\inf_{C_1(\delta_j)}I$.
 Pass to a further subsequence such that $\lim_{j\to\infty}z_j = z \in C$.
Then regardless of whether $z$ lies in the interior of $U$ or not, by (\ref{defi})
$\liminf_{j \to \infty} I(z_j)\geq I_d(z) \geq \inf_{C}I_d$, and we get the final upper bound

$$\limsup_{n\to \infty}\frac{1}{n}\log P_{0,\omega}^d\left(\frac{X_n}{n} \in C\right) \leq -\inf_{x \in C}I_d(x).$$

\medskip

\section{Proof of Theorem \ref{two} for the nearest neighbor case}
\label{s-four}
Here we consider the case in which the
jump range $R$ of the random walk $\{X_n:n\ge 0\}$ is nearest neighbor.
Define the {\it even lattice} as $\Z^d_{{\tiny even}}:=\{x \in \Z^d:
x_1+\ldots+x_d \, \mbox{is}\, \mbox{even}\}$. 
 Note that $\mathbb Z^d_{even}$ is
a free Abelian group which is isomorphic to $\mathbb Z^d$.
It therefore has a basis $f_1,\ldots f_d\in \mathbb Z^d_{even}$
and there is an isomorphism $h:\mathbb Z^d_{even}\to\mathbb Z^d$
such that $h(f_i)=e_i$ for $1\le i\le d$. It is obvious that
$h$ can be extended as an automorphism defined in $\mathbb R^d$.
Now, note that the random walk $\{Y_n:n\ge 0\}$ defined as

$$
Y_n:=h(X_{2n}),
$$ 
is a random walk in $\mathbb Z^d$ with finite, convex and symmetric jump range
$Q=h(R)$ and such that a neighborhood of the origin is contained in
its convex hull.
From Theorem \ref{two} for this class of random walks proved in
section 3, it follows that $\{Y_n:n\ge 0\}$ satisfies
a large deviation principle with a rate function $I$.
From this and the linearity of $h$ we conclude that the limit
\begin{equation}
\label{even preliminar function}
I_{even}(x):=I(h(x))=-\lim_{n\to\infty}\frac{1}{2n}\log\pi_{0,2n}(0,h^{-1}([2nh(x)])),
\end{equation}
exists $Q_\mu^d$-a.s, where $\pi_{n,m}(x,y)$ is the probability that the
random walk $\{X_n:n\ge 0\}$ jumps from time $n$ to time $m$
from site $x$ to site $y$. Furthermore, if $U:=\{x\in\mathbb R^d:|x|\le 1\}$, 
as in (\ref{defi}), one can define

\begin{equation}
\label{defi2}
I_{d,even}(x):=\left\{
\begin{array}{ll}
I_{even}(x) & {\rm for}\ x\in U^o\\
\liminf_{U^o\ni y\to x} I_{even}(y) & {\rm for}\ x\in \partial U \\
\infty & {\rm for}\ x\notin U,
\end{array}
\right.
\end{equation}
and $\{X_{2n}:n\ge 0\}$ satisfies a large deviation principle with
rate function $I_{even}$.

At this point, we need to extend the above large deviation principle
for the walk at even times, to all times taking into account the odd
number of steps of the random walk. The next lemma will be very useful for this objective. To do this, we first prove that
for each $x\in\mathbb R^d$ and each $g\in\mathcal{H}:=\left\{\sum_{i=1}^d c_i
x: c_i \in \{-1,0,1\}, x\in R\right\}$ we have that,

\begin{equation}
\label{odd}
I_{even}(x):=-\lim_{n\to\infty}\frac{1}{2n}\log\pi_{0,2n}(0,h^{-1}([2nh(x)])+g)
\qquad Q^d_\mu{\rm -a.s.}
\end{equation}
Note that to prove (\ref{odd}), it is enough to show that
for every  $g\in\mathcal{H}$ we have that,

\begin{equation}
\label{odd2}
\lim_{n\to\infty}\frac{1}{n}\log\tilde\pi_{0,n}(0,[nh(x)]+h(g))
=\lim_{n\to\infty}\frac{1}{n}\log\tilde\pi_{0,n}(0,[nh(x)]),
\end{equation}
where $\tilde \pi_{n,m}(x,y)$ is the probability that the random walk
$\{Y_n:n\ge 1\}$ jumps from time $n$ to time $m$ from site $x$ to site $y$.
The proof that the limit in the right-hand side of (\ref{odd2}) exists,
is a repetition of the proofs of Lemma \ref{s3lem} and
Proposition \ref{propos}, so we omit it. We just point out here  that in
 the proof of Lemma \ref{s3lem} we  need to replace the points $[nz]$, $[n_1x]$ and $[n_2x]$ 
by $[nz]+h$, $[n_1x]+h$ and $[n_2 x]+h$ respectively. On the
other hand, the equality in (\ref{odd2}) is established using
the uniform ellipticity of the walk and the Markov property.

Let us now see how to derive from (\ref{odd}) the large deviation principle for a random
walk with a nearest neighbor jump range $R$. Note that for any subset $A
\subseteq \R^d$ one has that

$$P_{0,\omega}\left(\frac{X_{2n+1}}{2n+1} \in A \right)=\sum_{i=1}^{2d} \pi_{0,1}(0,e_i)P_{e_i,\omega}\left(\frac{X_{2n}}{2n} \in A\right) = \sum_{i=1}^{2d} \pi_{0,1}(0,e_i)P_{0,\bar{\omega}}\left(\frac{X_{2n}}{2n} \in A-\frac{e_i}{2n}\right)$$
 where $\bar{\omega}=\{\omega_n: n \geq 1\}$ and $e_{i+d}=-e_i$ for $i=1, \ldots, d$. We will show that $P_{e_i,\omega}\left(\frac{X_{2n}}{2n} \in A\right)$ does not depend on $e_i$, regardless of whether $A$ is an open subset or a closed subset of $\R^d$ and we will use the result obtained in the even case. It is important to note that this argument can be used, even with $\bar{\omega}$, because the limit depends only on the distribution of $\omega$.

 Now, when $A=G$, where $G$ is an open subset of $\R^d$, we can 
follow the arguments used in the convex case, observing that for any
$x \in G$ and any $i \in \{1, \ldots, d\}$, $[nx]+e_i \in nG$, for 
$n$ large enough. On the other hand, if $A=C$, where $C$ is a compact subset
of $U_2^{\circ}$, note that

\begin{eqnarray*}
\limsup_{n \to \infty}\frac{1}{2n}\log P_{0, \bar{\omega}}\left(\frac{X_{2n}}{2n} \in C-\frac{e_i}{2n}\right)& \leq & \limsup_{n \to \infty}\sup_{x \in C-\frac{e_i}{2n}}\frac{1}{2n}\log \pi_{0,2n}(0,h^{-1}([2nh(x)])) \\
&= & \limsup_{n \to \infty}\sup_{x \in C}\frac{1}{2n}\log \pi_{0,2n}(0,h^{-1}
([2nh(x)-h(e_i)]))\\
&\leq& \limsup_{n \to \infty}\sup_{x \in C} \max_{g \in \mathcal{H}}\frac{1}{2n}\log \pi_{0,2n}(0,h^{-1}([2nh(x)])+g)\\
\end{eqnarray*}
However, by (\ref{odd}) the last expression is independent of $g$.

\medskip

\noindent {\bf Acknowledgments.} David Campos gratefully acknowledges the
support of the fellowship Consejo Nacional de Ciencia y Tecnolog\'\i a
number D-57080025,
Alexander Drewitz  of an ETH Fellowship,
Alejandro F. Ram\'\i rez
 of Fondo Nacional de Desarrollo Cient\'\i fico
y Tecnol\'ogico grant 1100298, Firas Rassoul-Agha
 of NSF Grant DMS-0747758 and Timo Sepp\"al\"ainen
  of NSF Grant DMS-100365. Also, the authors thank an anonymous
referee for several suggestions which led to an improved version of
this paper.


\begin{thebibliography}{}

\bibitem[ADHR10]{adhr10}
L. Avena, F. Redig, F. den Hollander.
{\it Large deviation principle for one-dimensional random walk in dynamic
  random environment: attractive spin-flips and simple symmetric exclusion.}
  Markov Process. Related Fields  {\bf 16},   139-168 (2010).
\bibitem[ADHR11]{adhr11}
L. Avena, F. Redig, F. den Hollander.
{\it Law of large numbers for a class of random walks in dynamic random
 environments.}  Electron. J. Probab.  {\bf 16}, 587-617 (2011).

\bibitem[ADSV11]{adsv11}
L. Avena, R. dos Santos, F. V\"ollering.
{\it Law of large numbers for a transient random walk driven by a symmetric
exclusion process.}
  arXiv:1102.1075 (2011).



\bibitem[BC06]{bc06}
G. Ben Arous, J. 
{\v{C}}ern{\'y}.
{\it Dynamics of trap models}.
 Mathematical statistical physics, 331-394,
 Elsevier B. V., Amsterdam, (2006).

\bibitem[DFGW89]{dfgw89}
A. De Masi, P.A. Ferrari, S. Goldstein, W.D. Wick.
{\it An invariance principle for
reversible Markov processes. Applications to random motions in random environments}. J.
Statist. Phys. {\bf 55}, 787-855 (1989).



\bibitem[DZ98]{dz98}
A. Dembo, O. Zeitouni.
{\it Large deviations techniques and applications}.
Springer-Verlag New York (1998).

\bibitem[DGRS12]{dgrs12}
A. Drewitz, J. G\"artner, A.F. Ram\'\i rez, R. Sun.
{\it Survival probability for a random walk on moving traps.}
In: J.D. Deuschel, B. Gentz, W. K\"onig, M.-K. von Renesse, M. Scheutzow,
U. Schmock (eds), Probability in Complex physical systems, Vol. 11,
pp. 119-158. Springer, Heidelberg (2012).


\bibitem[KL99]{kl99}
C.\ Kipnis, C.\ Landim.
{\it Scaling limits of interacting particle systems.}
 Springer-Verlag Berlin Heidelberg (1999).


\bibitem[L85]{l85}
T.\ Liggett.
{\it An improved subadditive ergodic theorem.}
 Ann.\ Probab. {\bf 13}, 1279--1285, (1985).


\bibitem[DHDSS11]{dhdss11}
F. den Hollander, R. dos Santos, V. Sidoravicius.
{\it Law of large numbers for non-elliptic random walks in dynamic random environments.}
arXiv:1103.2805 (2011).


\bibitem[KV08]{kv08}
E. Kosygina, S.R.S. Varadhan.
{\it Homogenization of Hamilton-Jacobi-Bellman equations with
 respect to time-space shifts in a
stationary ergodic medium.}
Comm. Pure Appl. Math. {\bf 61}, 816–847 (2008).

\bibitem[DMM86]{dmm86}
G. dal Maso, L. Modica.
{\it Nonlinear stochastic homogenization and ergodic theory.}
J. Reine Angew. Math. {\bf 368}, 28–42 (1986).


\bibitem[RSY11]{rsy11}
F. Rassoul-Agha, T. Sepp\"al\"ainen, A. Yilmaz.
{\it Quenched Free Energy and Large Deviations for Random Walks in Random Potentials}.
 arXiv:1104.3110 (2011).



\bibitem[R02]{r02}
F. Rezakhanlou.
{\it Continuum limit for some growth models}.
Stochastic Process. Appl. {\bf 101}, 1-41 (2002).

\bibitem[RT00]{rt00}
F. Rezakhanlou, J. E. Tarver.
{\it Homogenization for stochastic Hamilton-Jacobi equations}.
 Ration. Mech. Anal. {\bf 151}, 277-309 (2000).



\bibitem[R97]{r97}
R. T. Rockafellar.
{\it Convex analysis}.
Princeton University Press  (1997).



\bibitem[R06]{r06}
J. Rosenbluth.
{\it  Quenched large deviations for multidimensional random walk in random
  environment: a variational formula.}
Thesis (Ph.D.)–New York University (2006).


\bibitem[S99]{s99}
T. Sepp\"al\"ainen.
{\it Existence of hydrodynamics for the totally asymmetric simple
  $K$-exclusion process.}
  Ann. Probab.  {\bf 27},   361-415 (1999).

\bibitem[So99]{so99}
P. Souganidis.
{\it Stochastic homogenization of Hamilton-Jacobi equations and
 some applications.}
Asymptot. Anal. {\bf 20}, 1–11 (1999).


\bibitem[S98]{s98}
A.S. Sznitman.
{\it Brownian motion, Obstacles and Random Media.}
  Springer-Verlag Berlin Heidelberg (1998).


\bibitem[S04]{s04}
A.S. Sznitman.
{\it Topics in random walks in random environment.}
School and Conference on Probability Theory,  203-266,
ICTP Lect. Notes, XVII, Abdus Salam Int. Cent. Theoret. Phys., Trieste, (2004).

\bibitem[V03]{v03}
S.\ R.\ S. Varadhan. {\it Large deviations for random walks in a random
  environment}, Comm. Pure Appl. Math. {\bf 56}, 1222-1245
(2003).


\bibitem[Y08]{y08}
A. Yilmaz.
{\it Large deviations for random walk in a random environment.}
 Thesis (Ph.D.)–New York University (2008).

\bibitem[Y09]{y09}
A. Yilmaz.
{\it Large deviations for random walk in a space-time product environment.}
  Ann. Probab.  {\bf 37}, 189-205 (2009).


\bibitem[Y09-2]{y09-2}
A. Yilmaz.
{\it Quenched large deviations for random walk in a random environment.}
Comm. Pure Appl. Math. {\bf 62}, 1033–1075 (2009).

\bibitem[Z06]{z06}
O. Zeitouni.
{\it Random walks in random environments.}
J. Phys. A {\bf 39}, R433–R464 (2006).


\bibitem[Z98]{z98}
M. Zerner.
{\it Lyapounov exponents and quenched large deviations
 for multidimensional random walk in random environment.}
Ann. Probab. {\bf 26},  1446–1476 (1998).

\end{thebibliography}
\end{document}